\documentclass[12pt]{amsart}
\usepackage{mathrsfs}
\usepackage{amsmath,amssymb,amsfonts,latexsym,txfonts}
\usepackage{eucal}
\usepackage{fancyhdr}

\newlength{\originalbase}
\newtheorem{theorem}{Theorem}[section]
\newtheorem{proposition}[theorem]{Proposition}
\newtheorem{lemma}[theorem]{Lemma}

\newtheorem{corollary}[theorem]{Corollary}

\input amssym.def
 
\input amssym.tex

\renewcommand{\SS}{\mathbb{S}}
\newcommand{\RR}{\mathbb{R}}

\newcommand{\LL}{\mathcal{L}}
\newcommand{\KK}{\mathscr{K}}

\newcommand{\FF}{\mathscr{F}}

\begin{document}

\title{Containment and inscribed simplices}

\author{Daniel A. Klain}
\address{Department of Mathematical Sciences,
University of Massachusetts Lowell,
Lowell, MA 01854 USA
}

\email{Daniel\_{}Klain@uml.edu}

\begin{abstract}
Let $K$ and $L$ be compact convex sets in $\RR^n$.  
The following two statements are shown to be equivalent:
\begin{enumerate}
\item[\bf (i)] For every polytope $Q \subseteq K$ having at most $n+1$ vertices, 
$L$ contains a translate of $Q$.
\item[\bf (ii)] $L$ contains a translate of $K$.
\end{enumerate}
Let $1 \leq d \leq n-1$.  
It is also shown that the following two statements are equivalent:
\begin{enumerate}
\item[\bf (i)] For every polytope $Q \subseteq K$ having at most $d+1$ vertices, 
$L$ contains a translate of $Q$.
\item[\bf (ii)] For every $d$-dimensional subspace $\xi$, the orthogonal projection 
$L_\xi$ of the set $L$ contains a translate of the corresponding projection $K_\xi$ of the set $K$.
\end{enumerate}
It is then shown that,
if $K$ is a compact convex set in $\RR^n$ having at least $d+2$ exposed points,
then there exists a compact convex set $L$ such that every $d$-dimensional orthogonal projection
$L_\xi$ contains a translate of
the projection $K_\xi$, while $L$ does not contain a translate of $K$. 
In particular, if $\dim K > d$, then there exists $L$ 
such that every $d$-dimensional projection $L_\xi$ contains a translate of
the projection $K_\xi$, while $L$ does not contain a translate of $K$.
\end{abstract}

\maketitle

This note addresses questions related to following general problem:
Consider two compact convex subsets $K$ and $L$ of 
$n$-dimensional Euclidean space.  
Suppose that, for a given dimension $1 \leq d < n$,
every $d$-dimensional orthogonal projection (shadow) of $L$
contains a translate of the corresponding projection of $K$.
Under what conditions does it follow that the original set $L$ contains a translate of $K$?  
In other words, if $K$ can be translated to ``hide behind" $L$ from any perspective,
does it follow that $K$ can ``hide inside" $L$?  

This question is easily answered when a sufficient degree of symmetry is imposed.
For example, a support function argument implies that the answer is {\em Yes}
if {\em both} of the bodies $K$ and $L$ are centrally symmetric.
It is also not difficult to show that if every $d$-projection of $K$ (for some $1 \leq d < n$)
can be translated into the corresponding shadow of an orthogonal $n$-dimensional box $C$,
then $K$ fits inside $C$ by some translation, since one needs only to check that the widths
are compatible in the $n$ edge directions of $C$.  A similar observation applies if $C$
is a parallelotope (an affine image of a box), a cylinder (the product of
an $(n-1)$-dimensional compact convex set with a line segment), or a similarly decomposable
product set; see also \cite{Klain-Shadow}).

For more general classes of convex bodies the situation is quite different.  
Given any $n > 1$ and $1 \leq d \leq n-1$,
it is possible to find convex bodies $K$ and $L$ in $\RR^n$ 
such every $d$-dimensional orthogonal projection (shadow) of $L$
contains a translate of the corresponding projection of $K$,
even though $K$ has {\em greater volume} than $L$ (and so certainly could not fit inside $L$).
For a detailed example of this volume phenomenon, see \cite{Klain-Shadow}.

In \cite{Lut-contain} Lutwak uses Helly's theorem to
prove that, if every $n$-simplex containing $L$ also contains a translate of $K$, then
$L$ contains a translate of $K$.  In the present note we describe a dual result, by which
the question of containment is related to properties of the {\em inscribed} simplices (and more
general polytopes) of the bodies $K$ and $L$.  We then generalize these containment (covering) theorems 
in order to reduce questions about shadow (projection) covering to questions about inscribed simplices and
related polytopes.  Specifically we establish the following:
\begin{enumerate}
\item[(1)]  Let $K$ and $L$ be compact convex sets in $\RR^n$.  
The following are equivalent:
\begin{enumerate}
\item[\bf (i)] For every polytope $Q \subseteq K$ having at most $n+1$ vertices, 
$L$ contains a translate of $Q$.
\item[\bf (ii)] $L$ contains a translate of $K$.
\end{enumerate}
(Theorem~\ref{inscr})\\

\item[(2)] Let $K$ and $L$ be compact convex sets in $\RR^n$, and let $1 \leq d \leq n-1$.  
The following are equivalent:
\begin{enumerate}
\item[\bf (i)] For every polytope $Q \subseteq K$ having at most $d+1$ vertices, 
$L$ contains a translate of $Q$.
\item[\bf (ii)] For every $d$-dimensional subspace $\xi$, the orthogonal projection 
$L_\xi$ contains a translate of $K_\xi$.
\end{enumerate}
(Theorem~\ref{i-shadow})\\

\item[(3)] Let $1 \leq d \leq n-1$.
If $K$ is a compact convex set in $\RR^n$ having at least $d+2$ exposed points,
then there exists a compact convex set $L$ such that every $d$-dimensional orthogonal projection
$L_\xi$ contains a translate of
the projection $K_\xi$, while $L$ does not contain a translate of $K$ itself. 
(Theorem~\ref{toomany-d})
\end{enumerate}
In particular, if $\dim K > d$, then there exists $L$ 
such that every $d$-shadow $L_\xi$ contains a translate of
the shadow $K_\xi$, while $L$ does not contain a translate of $K$.

In this note we address the existence of a compact convex set $L$, 
whose shadows can cover those of a given set $K$, without containing a translate of $K$ itself.
A reverse question is addressed in \cite{Klain-Circ}:  Given a body $L$, does there necessarily exist $K$ so that
the shadows of $L$ can cover those of $K$, while $L$ does not contain a translate of $K$?
These containment and covering problems are special cases of
the following more general question: Under what conditions 
will a compact convex set necessarily contain a translate or otherwise congruent copy of another?  
Progress on different aspects of this general question also appears in the work of
Gardner and Vol\v{c}i\v{c} \cite{Gard-Vol},
Groemer \cite{Gro-proj},
Hadwiger \cite{Had3,Had4,Had-proj,Lincee,Santa}, 
Jung \cite{Bonn2,Webster}, 
Lutwak \cite{Lut-contain}, 
Rogers \cite{Rogers}, 
Soltan \cite{Soltan},
Steinhagen \cite[p. 86]{Bonn2}, 
Zhou \cite{Zhou1,Zhou2}, and many others (see also \cite{Gard2006,Klain-Circ,Klain-Shadow}).

\section{Background} 
Denote $n$-dimensional Euclidean space by $\RR^n$, and let $\SS^{n-1}$ denote the set of
unit vectors in $\RR^n$; that is, the unit $(n-1)$-sphere centered at the origin.

Let $\KK_n$ denote the set of compact convex subsets of $\RR^n$.  If $u$ is a unit vector in
$\RR^n$, denote by $K_u$ the orthogonal projection of a set $K$ onto the subspace $u^\perp$.
More generally, if $\xi$ is a $d$-dimensional subspace of $\RR^n$, denote by $K_\xi$
the orthogonal projection of a set $K$ onto the subspace $\xi$.  The boundary of a compact
convex set $K$ will be denoted by $\partial K$.

Let $h_K: \RR^n \rightarrow \RR$ denote the support function of a compact convex set $K$;
that is,
$$h_K(v) = \max_{x \in K} x \cdot v$$
For $K, L \in \KK_n$, we have $K \subseteq L$ if and only if $h_K \leq h_L$.  
If $\xi$ is a subspace of $\RR^n$ then the support function $h_{K_\xi}$ is given by the restriction
of $h_K$ to $\xi$ (see also \cite[p. 38]{red}).

If $u$ is a unit vector in
$\RR^n$, denote by $K^u$ the support set of $K$ in the direction of $u$; that is,
$$K^u = \{x \in K \; | \; x \cdot u = h_K(u) \}.$$
If $P$ is a convex polytope, then $P^u$ is the face of $P$ having $u$ in its outer normal cone.
A point $x \in \partial K$ is an {\em exposed point} of $K$ if $x = K^u$ for some direction $u$.
In this case, the direction $u$ is said to be a {\em regular unit normal} to $K$.
If $K$ has non-empty interior, then the regular unit normals to $K$ are dense in the 
unit sphere $\SS^{n-1}$ (see \cite[p. 77]{red}).

Suppose that $\mathscr{F}$ is a family of compact convex sets in $\RR^n$.
Helly's Theorem \cite{Bonn2,Lincee,red,Webster} asserts that, if every $n+1$ sets
in $\mathscr{F}$ share a common point, then the entire family shares a common point.
In \cite{Lut-contain} Lutwak used Helly's theorem to
prove the following fundamental criterion for whether a set $L \in \KK_n$ contains 
a translate of another compact convex set $K$.   
\begin{theorem}[Lutwak's Containment Theorem] 
Let $K,L \in \KK^n$.  
The following are equivalent:
\begin{enumerate}
\item[\bf (i)] For every simplex $\Delta$ such that $L \subseteq \Delta,$ 
there exists $v \in \RR^n$ such that 
$K + v \subseteq \Delta$.
\item[\bf (ii)] There exists $v_0 \in \RR^n$ such that
$K + v_0 \subseteq L$.
\end{enumerate}
\label{lut}
\end{theorem}
In other words, if every $n$-simplex containing $L$ also contains a translate of $K$, then
$L$ contains a translate of $K$.  

\section{Inscribed polytopes and shadows}

The following theorem provides an {\em inscribed} polytope counterpart to Lutwak's theorem.
\begin{theorem}[Inscribed Polytope Containment Theorem]
Let $K,L \in \KK^n$.  The following are equivalent:
\begin{enumerate}
\item[\bf (i)] For every polytope $Q \subseteq K$ having at most $n+1$ vertices, 
there exists $v \in \RR^n$ such that 
$Q + v \subseteq L$.
\item[\bf (ii)] There exists $v_0 \in \RR^n$ such that
$K + v_0 \subseteq L$.
\end{enumerate}
\label{inscr}
\end{theorem}

\begin{proof}
The implication {\bf (ii)} $\Rightarrow$ {\bf (i)} is
obvious.  We show that {\bf (i)} $\Rightarrow$ {\bf (ii)}.

Note that $x + v \in L$ if and only if $v \in L-x$.
If $x_0, x_1, \ldots, x_{n} \in K$, let $Q$ denote the convex hull of these points.
Note that $Q$ has at most $n+1$ vertices.  By
the assumption {\bf (i)} there exists $v$ such that
$Q + v \subseteq L$.  In other words, $x_i + v \in L$ for each $i$, so that
\begin{equation}
v \in \bigcap_{i=0}^n (L-x_i).
\label{found}
\end{equation}

Let $\FF = \{L-x \; | \; x \in K\}$.  By~(\ref{found}), $\FF$ is a family of compact convex sets
that satisfies the intersection condition of Helly's theorem \cite{red,Webster}.  
Hence there exists
a point $v_0$ such that
$$v_0 \in \bigcap_{x \in K} (L-x).$$
In other words, $x+v_0 \in L$ for all $x \in K$, so that $K+v_0 \subseteq L$.
\end{proof}

\begin{corollary}
Suppose that $K, L \in \KK_n$ have {\em non-empty interiors}.
If every {\em simplex} contained in $K$ can be translated inside
$L$, then $K$  can be translated inside $L$.
\label{interior}
\end{corollary}

\begin{proof}  The proof is the same as that of Theorem~\ref{inscr}, except
that we must address the case in which the points 
$x_0, x_1, \ldots, x_{n} \in K$ are affinely dependent (and are not the vertices
of a simplex).

In this case, since $K$ has interior, perturbations
of these points by a small distance $\epsilon > 0$
will yield the vertices of a simplex and a vector $v_\epsilon$ such that~(\ref{found}) holds
for the perturbed points.
As $\epsilon \rightarrow 0$ a vector $v$ is obtained so that~(\ref{found}) holds for the
original points $x_0, x_1, \ldots, x_{n}$ as well, since $L$ is compact.  Helly's theorem
now applies, as in the previous proof.
\end{proof}

Theorem~\ref{inscr} is now generalized to address covering of lower-dimensional shadows.
\begin{theorem}[Generalized Inscribed Polytope Containment Theorem] \hfill\\
Let $K,L \in \KK^n$, and suppose $1 \leq d \leq n$.
The following are equivalent:
\begin{enumerate}
\item[\bf (i)] For every polytope $Q \subseteq K$ having at most $d+1$ vertices, 
there exists $v \in \RR^n$ such that 
$Q + v \subseteq L$.
\item[\bf (ii)] For every $d$-dimensional subspace $\xi$, there exists $v \in \xi$
such that $K_\xi + v \subseteq L_\xi$.
\end{enumerate}
\label{i-shadow}
\end{theorem}
When $K$ and $L$ have non-empty interiors, this theorem can be reformulated in
the following way:
if every $d$-simplex contained in $K$ can be translated into
$L$, then every $d$-shadow of $K$ can be translated into the corresponding $d$-shadow of $L$,
and vice versa.  
In this case a perturbation argument applies, as in the proof of Corollary~\ref{interior}.

The next three lemmas will be used to prove
Theorem~\ref{i-shadow}.

\begin{lemma}  
Let $T$ be an $n$-simplex, and let $Q$ be a polytope in $\RR^n$ having at most $n$ vertices.
Suppose that, for every unit vector $u$, there exists $v \in u^\perp$
such that $Q_u + v \subseteq T_u$.  Then there exists $v_0 \in \RR^n$ such that
$Q + v_0 \subseteq T$.
\label{spx}
\end{lemma}

\begin{proof}  Since $T$ has interior, $\epsilon Q$ can be translated 
inside $T$ for sufficiently small $\epsilon > 0$.  Let $\hat{\epsilon}$ 
denote the maximum of all such $\epsilon > 0$.  We will show that
$\hat{\epsilon} \geq 1$, thereby proving the lemma.

Without loss of generality, translate $T$ so that $\hat{\epsilon}Q \subseteq T$.  If 
$\hat{\epsilon}Q$ does not intersect a given facet $F$ of $T$, 
then some translate of $\hat{\epsilon}Q$ lies in the {\em interior} of $T$.
This violates the maximality of $\hat{\epsilon}$.  It follows that $\hat{\epsilon}Q$ must meet
every facet of $T$.  In particular, the vertex set of $\hat{\epsilon}Q$ must meet every facet of $T$.
Since $\hat{\epsilon}Q$ has at most $n$ vertices, while $T$ has $n+1$ facets, some vertex of $\hat{\epsilon}Q$
must meet a face $\sigma$ of $T$ having co-dimension 2, where $\sigma = F_1 \cap F_2$, the intersection
of two facets of $T$.

Let $\ell$ denote the line segment (i.e.~the edge) complementary to $\sigma$ in the boundary $\partial T$ (so that
$T$ is the convex hull of the union $\ell \cup \sigma$).  If $v \in \RR^n$ points in the direction of $\ell$,
then $T_v$ is an $(n-1)$-simplex.  Moreover, every facet of $T_v$ except one is exactly the projection of
a facet of $T$, while $(F_1)_v = (F_2)_v = T_v$.  The remaining facet of $T_v$ is the projection
$\sigma_v$ of the ridge $\sigma$ in $T$.  Since $\hat{\epsilon}Q$ meets every facet of $T$, as well as
the ridge $\sigma$, the projection $\hat{\epsilon}Q_v$ meets every facet of $T_v$, 
and is therefore inscribed (maximally) in $T_v$.  Therefore, if $\epsilon > \hat{\epsilon}$, then $\epsilon Q_v$
cannot be translated inside $T_v$.  Since every shadow $Q_v = 1 Q_v$ of $Q$ can be translated inside 
the corresponding shadow of $T_v$ (by hypothesis), it follows that $\hat{\epsilon} \geq 1$.
\end{proof}

\begin{lemma}  
Let $L \in \KK_n$, and let $Q$ be a polytope in $\RR^n$ having at most $n$ vertices.
If every shadow $L_u$ contains a translate of the corresponding shadow $Q_u$,
then $L$ contains a translate of $Q$.
\label{spx1}
\end{lemma}

\begin{proof}  
Let $T$ be an $n$-simplex that contains $L$.
Since $Q_u$ can be translated inside the corresponding shadow $L_u$, for each $u$, it follows
that $Q_u$ can be translated inside the corresponding shadow $T_u \supseteq L_u$ as well.  
By Lemma~\ref{spx}, $Q$ can be translated inside $T$.  
Since this holds for every $n$-simplex $T \supseteq L$, 
Lutwak's Theorem~\ref{lut} implies that $L$ contains a translate of $Q$.
\end{proof}

\begin{lemma}  
Let $L \in \KK_n$, and let $Q$ be a polytope in $\RR^n$ having at most $d+1$ vertices, where $d < n$.
Suppose that, for every $d$-dimensional subspace $\xi$, there exists $v \in \xi$
such that $Q_\xi + v \subseteq L_\xi$.  Then there exists $v_0 \in \RR^n$ such that
$Q+ v_0 \subseteq L$.
\label{spx2}
\end{lemma}

\begin{proof} Fix $d$ and proceed by induction on $n$, starting with the case $n = d+1$,
which follows from Lemma~\ref{spx1}.  

Now suppose that Lemma~\ref{spx2} is true for $n \leq d+i$.   If $n = d+i+1$, then 
each projection $Q_u$ also has at most $d+1$ vertices, 
The induction assumption (in the lower dimensional space $u^\perp$)
applies to $Q_u$, so that $Q_u$ can be translated inside $L_u$ for all $u$.
Because $Q$ has at most $d+1 \leq n$ vertices, 
Lemma~\ref{spx1} implies that $Q$ can be translated inside $L$.
\end{proof}

We now prove Theorem~\ref{i-shadow}.
\begin{proof}[{\bf\em Proof of Theorem~\ref{i-shadow}}]  To begin suppose that {\bf (i)} holds.
If $Q \subseteq K_\xi$ has at most $d+1$ vertices, then $Q$ is the projection of a polytope
$\tilde{Q} \subseteq K$ having at most $d+1$ vertices.  By {\bf (i)} there exists $v \in \RR^n$ such that
$\tilde{Q} + v \subseteq L$.  By the linearity of orthogonal projection
it follows that $Q + v_{\xi} \subseteq L_\xi$.  The assertion {\bf (ii)} now follows
from Theorem~\ref{inscr} applied inside the subspace $\xi$.

To prove the converse, suppose that {\bf (ii)} holds.  Let $Q \subseteq K$ be a polytope
with at most $d+1$ vertices. 
For each $\xi$ there exists $w \in \xi$ such that $K_\xi + w \subseteq L_\xi$, by {\bf (ii)}.
Since $Q \subseteq K$, we have $Q_\xi + w \subseteq K_\xi + w \subseteq L_\xi$ as well.
It follows from Lemma~\ref{spx2} that
there exists  
$v \in \RR^n$ such that $Q + v \subseteq L$.  
\end{proof}

Webster \cite[p. 301]{Webster} shows that if every triangle inside a compact convex set $K$ can be translated inside
a compact convex set $L$ of the {\em same diameter} as $K$, then $K$ can itself be translated inside $L$.  
Combining this observation with Theorem~\ref{i-shadow} yields the following corollary.
\begin{corollary}  Let $K, L \in \KK_n$, and let $d \geq 2$.  Suppose that 
every $d$-dimensional shadow $L_\xi$ contains a translate of the corresponding shadow $K_\xi$.
If $K$ and $L$ have the same diameter, then $L$ contains a translate of $K$.
\end{corollary}

Webster's observation can be generalized in other ways via Theorem~\ref{i-shadow}.
Denote by $W(K)$ the {\em mean width} of the body $K$, taken over all directions in $\RR^n$.
If $h_K$ is the support function of $K$, then 
$$W(K) = \frac{2}{n \omega_n} \int_{\SS^{n-1}} h_K(u) \; du,$$
where $\omega_n$ is the volume of the $n$-dimensional Euclidean unit ball.  Evidently $W(K)$ is
strictly monotonic, in the sense that $W(K) \leq W(L)$ whenever $K \subseteq L$, with
equality if and only if $K = L$.  (This follows from the fact that a compact convex set is uniquely determined
by its support function \cite{Bonn2,red}.)  An alternative way to compute the mean width is
given by the following Kubota-type formula \cite{Bonn2,Lincee,red}:
\begin{equation}
W(K) = \int_{G(n,2)} W(K_\xi) \; d\xi,
\label{kub}
\end{equation}
where $G(n,2)$ is the Grassmannian of 2-dimensional subspaces of $\RR^n$, and the integral
is taken with respect to Haar probability measure. 
\begin{corollary}  
Let $K, L \in \KK_n$.
Suppose that every triangle inside $K$ can be translated inside $L$. 
If $K$ and $L$ have the same mean width, then 
$K$ and $L$ are translates.  
\end{corollary}

\begin{proof} By Theorem~\ref{i-shadow}, every 2-dimensional shadow of $K$ can be translated inside
the corresponding shadow of $L$.  It follows that $W(K_\xi) \leq W(L_\xi)$ for each 2-subspace $\xi$.
If $W(K) = W(L)$, then~(\ref{kub}) and the monotonicity of $W$ yields
$$W(K) = \int_{G(n,2)} W(K_\xi) \; d\xi
\leq \int_{G(n,2)} W(L_\xi) \; d\xi = W(L) = W(K),$$
so that equality $W(K_\xi) = W(L_\xi)$ holds in every 2-subspace $\xi$.  The strictness of monotonicity
for $W$ now implies that each $L_\xi$ is a {\em translate} of $K_\xi$.  

A well-known theorem asserts that if $K$ and $L$ have 
translation-congruent 2-dimensional projections, then
$K$ and $L$ are translates
(see, for example, \cite[p. 100]{Gard2006} or \cite{Gro-proj,Had-proj,Rogers}).
\end{proof}

The concept of mean width can be generalized to quermassintegrals (mean $d$-volumes of $d$-dimensional
shadows).  The previous argument (combining Theorem~\ref{i-shadow} with 
monotonicity, Kubota formulas, and the homothetic projection theorem)
generalizes to give the following. 
\begin{corollary} Let $K, L \in \KK_n$, and let $d \geq 2$.
Suppose that every $d$-simplex inside $K$ can be translated inside $L$.  
If $K$ has the same $m$-quermassintegral as $L$, 
for some $1 \leq m \leq d$, then 
$K$ and $L$ are translates.  
\end{corollary}

The previous corollary does not hold for $m > d$.  For example, there exist convex bodies
$K$ and $L$ in $\RR^3$ such that $L$ contains a translate of every triangle inside $K$, even though $L$ has {\em strictly smaller} volume than $K$.  Explicit examples of this phenomenon are described in \cite{Klain-Shadow}.
In this case every $2$-shadow $K_\xi$ can be translated inside the corresponding shadow $L_\xi$ (by Theorem~\ref{i-shadow}), while the (Euclidean) volumes of $L$ and $K$ satisfy $V(L) < V(K)$.  
This implies that $K$ and $L$ are not homothetic.  Now dilate $L$ sufficiently 
so that $V(L) = V(K)$.  The triangle covering condition is preserved, but $K$ and $L$ are not translates.

\section{Most objects may be hidden without being covered}

We have shown that, if the $d$-shadows of
a compact convex set $L$ cover the $d$-shadows of a polytope $Q$
having at most $d+1$ vertices, then $L$ contains a translate of $Q$.  What if $Q$ has more vertices?
What if $Q$ is replaced by a more general compact convex set $K$?  It turns out that adding one additional vertex
changes the story.

Consider, for example, a regular tetrahedron $\Delta$ in $\RR^3$.  Let $Q$ be a 
planar quadrilateral with one vertex from the relative interior of each facet of $\Delta$.  
Since $Q$ does not meet any edge of $\Delta$, 
every 2-shadow of $Q$ has a translate inside 
the {\em interior} of the corresponding 2-shadow of $\Delta$.  
By a standard compactness argument, 
there is an $\epsilon > 1$ such that every 2-shadow of $\epsilon Q$ can be translated inside 
the corresponding 2-shadow of $\Delta$.  But $Q$ already meets every facet of $\Delta$, so the simplex
$\Delta$ cannot contain any translate of $\epsilon Q$.

More generally, we will show that if $K \in \KK_n$ has more than $d+1$ exposed points, 
then there exists $L \in \KK_n$ whose $d$-shadows contain translates of
the corresponding $d$-shadows of $K$,
while $L$ does not contain a translate of $K$.

\begin{lemma} Let $\Delta$ be an $n$-simplex, and let $K \subseteq \Delta$ be a compact convex set.
Suppose that $K \cap F = \emptyset$ for every face $F$ of $\Delta$ such that $\dim(F) \leq n-2$.
Then
\begin{enumerate}
\item[\bf (i)] For each $u \in \SS^{n-1}$, the 
projection $K_u$ can be translated inside the interior of $\Delta_u$.\\

\item[\bf (ii)] There exists $\epsilon > 1$ such that, for each $u$, 
the projection $\Delta_u$ contains a translate of
$\epsilon K_u$.
\end{enumerate}
\label{um}
\end{lemma}
Note that the value $\epsilon$ in {\bf (ii)} is independent of the direction $u$.

\begin{proof} 

Since $K \subseteq \Delta$, each $K_u \subseteq \Delta_u$.  Suppose that some projection $K_u$ cannot 
be translated into
the interior of $\Delta_u$.  In this case, $K_u$ meets the boundary $\partial \Delta_u$ in supporting directions
$u_0, \ldots, u_k \in u^\perp \cap \SS^{n-1}$ such that the origin $o$ lies
in the relative interior of the convex hull of $u_0, \ldots u_k$; that is,
\begin{align}
a_0 u_0 + \cdots + a_k u_k = o,
\label{subsimp}
\end{align}
where each $a_i > 0$ and $a_0 + \cdots + a_k = 1$.  Moreover, by
Caratheodory's Theorem, applied in the $(n-1)$-dimensional space $u^\perp$,
we can assume that $k \leq n-1$. 
This means that 
$$h_K(u_i) = h_{K_u}(u_i) = h_{\Delta_u}(u_i) = h_\Delta(u_i),$$
for each $u_i$.  Because $k < n$, no $k+1$ facet normals of an $n$-simplex $\Delta$
can satisfy~(\ref{subsimp}).  Therefore, at least one of the directions $u_i$ is not a facet normal
of $\Delta$, so that $K$ must meet an $(n-2)$-dimensional face of $\Delta$, contradicting
the hypothesis of the lemma.



This proves {\bf (i)}.  

Since the interior of each $\Delta_u$ contains a translate of $K_u$,
there exists $\epsilon_u > 1$ such that $\epsilon_u K$ can be translated inside $\Delta_u$.
Let $\epsilon = \inf_u \epsilon_u$, and let $\{u_i\}$ be a sequence of unit vectors
such that $\epsilon_i = \epsilon_{u_i}$ converge to $\epsilon$.  Since the unit sphere is
compact, we can pass to a subsequence as needed, and assume without loss of generality that
$u_i \rightarrow v$ for some unit vector $v$.

Since $\epsilon_v > 1$, we can translate $K$ and $\Delta$ so that $o \in K_v \subseteq \Delta_v$,
where the origin $o$ now lies in the interior of $\Delta$.
If $\alpha = \frac{1+ \epsilon_v}{2}$, then $\alpha K_v$ lies in the relative interior of $\Delta_v$, so that
their support functions satisfy $\alpha h_K(x) < h_{\Delta}(x)$ for all unit vectors $x \in v^\perp$.  Since
support functions are uniformly continuous on the unit sphere, and since $u_i \rightarrow v$, we have
$\alpha h_K(x) < h_{\Delta}(x)$ for all $x \in u_i^\perp$ for $i$ sufficiently large.  This means that
$\alpha K_{u_i}$ lies in the relative interior of $\Delta_{u_i}$ for large $i$, so that
$\alpha < \epsilon_i$ as well.  Taking limits, we have $1 < \alpha \leq \epsilon$.
Since $\epsilon > 1$,
the assertion {\bf (ii)} now follows.
\end{proof}

A set $C \subseteq \SS^{n-1}$ is a closed spherical convex set 
if $C$ is an intersection of closed hemispheres.
The polar dual $C^*$ is defined by
$$C^* = \{u \in \SS^{n-1} \; | \; u \cdot v \leq 0 \hbox{ for all } v \in C \}.$$
If $x \in C \cap C^*$ then $x \cdot x = 0$.  This is impossible for a unit vector $x$, so we 
have $C \cap C^* = \emptyset$.  Recall also that $C^{**} = C$.  See, for example, \cite{red,Webster}.  
(Note that one can identify $C$ with the cone
obtained by taking all nonnegative linear combinations in $\RR^n$ of points in $C$, 
taking the polar dual in this context, and then intersecting with the sphere once again.)

\begin{lemma} Let $C$ be a closed spherical convex set in $\SS^{n-1}$.  
Then there exists a unit vector $v \in -C \cap C^*$.  

Moreover, if $C$ has dimension $j \geq 0$ and lies in
the interior of a hemisphere, then $-C \cap C^*$ also has dimension $j$.
\end{lemma}

\begin{proof} Since $C \cap C^* = \emptyset$,
there is a hyperplane $H = v^\perp$ through the origin in $\RR^n$ that separates them.
Let $H^+$ and $H^-$ denote the closed hemispheres bounded by $H \cap \SS^{n-1}$, \
labelled so that $v \in H^+$, and so that
$C \subseteq H^-$ and $C^* \subseteq H^+$.

Since $C \subseteq H^- \subseteq \{v\}^*$, we have $v \in C^*$.  (Polar duality reverses inclusion relations.)
Meanwhile, $C^* \subseteq H^+ =  -\{v\}^* = \{-v\}^*$, so that $-v \in C^{**} = C$, and $v \in -C$. 
Conversely, if $v \in -C \cap C^*$ then $v^\perp$ separates $C$ and $C^*$.

If $C$ has dimension $j \geq 0$ 
and lies in the interior of a hemisphere, then $C^*$ has interior, and the set $C^* \cap -C$ 
consists of  all $v$ such
that $v^\perp$ separates $C$ and $C^*$, a set of dimension $j$ as well. 
\end{proof}

\begin{theorem}  If $K \in \KK_n$ has dimension $n$, 
then there exist regular unit normal vectors $u_0, \ldots, u_n$,
at distinct exposed points $x_0, \ldots, x_n$ on the boundary of $K$,
such that $u_0, \ldots, u_n$ are the outward unit normals vectors of
some $n$-dimensional simplex in $\RR^n$.
\label{got-normals}
\end{theorem}
Note that Theorem~\ref{got-normals} is trivial if $K$ is smooth and strictly convex, where each supporting
hyperplane of $K$ meets $K$ at a single boundary point, and each boundary point has exactly one
supporting hyperplane.  In this case, {\em any} circumscribing $n$-simplex for $K$ will do.

If $K$ is a polytope, then Theorem~\ref{got-normals} is again easy to prove, since each exposed
point (vertex) of $K$ has a unit outward normal cone with interior in the unit sphere, and these
interiors fill the sphere except for a set of measure zero.  Once again we can take any circumscribing simplex $S$
for $K$,
and then make small perturbations of each facet normal so the each facet of $S$ meets a different vertex of $K$.

The following more technical argument verifies Theorem~\ref{got-normals} for arbitrary $K \in \KK_n$ having
dimension $n$ (i.e. having non-empty interior).
\begin{proof}[{\bf\em Proof of Theorem~\ref{got-normals}}]  
If $x$ lies on the boundary of $K$, denote by $N(K,x)$ the outward unit normal cone
to $K$ at $x$; that is,
$$N(K,x) = \{u \in \SS^{n-1} \; | \; x \cdot u = h_K(u) \}.$$

Let $u_0$ be a regular unit normal at the exposed point $x_0 = K^{u_0}$.   
By the previous lemma,
we can choose $u_0$ in the normal cone 
$N_0 = N(K,x_0)$ so that $u_0 \in N(K,x_0) \cap -N(K, x_0)^*$.

Since $K$ has dimension $n$, the normal cone $N_0$ lies in an open hemisphere.  
Recall that regular unit normal vectors to $K$ are dense in the unit sphere $\SS^{n-1}$ (see \cite[p. 77]{red}).
It follows that we can choose $u_1, x_1, N_1$ similarly, so that $u_1$ lies outside $N_0$
and so that $\{u_0, u_1\}$ are linearly independent.  Once again $N_1$ lies inside an open
hemisphere. 

Having chosen $u_i, x_i, N_i$ in this manner, for $i = 0, \ldots, k$, where $k < n-1$, the union
$N_0 \cup \ldots \cup N_k$ cannot cover the sphere, because each is a closed subset of an open hemisphere,
and the $S^{n-1}$ is not the union of $n-1$ open hemispheres.  It follows that
$$X = \SS^{n-1} - (N_0 \cup \ldots \cup N_k)$$
is a nonempty open subset of $\SS^{n-1}$.  Since regular unit normals to $K$ are dense in the sphere,
we can choose $u_{k+1} \in X$ so that $x_{k+1}$ is disjoint from the previous choices of $x_i$,
and such that $u_0, \ldots, u_{k+1}$ are linearly independent.

Continuing in this manner, we obtain a linearly independent set $u_0, \ldots, u_{n-1}$ of regular unit
normals at distinct exposed points $x_0, \ldots, x_{n-1}$ of $K$.  Since the unit normals 
$u_0, \ldots, u_{n-1}$ are independent,
the origin $o$ does {\em not} lie in their convex hull.  Therefore, there exists an open hemisphere containing
$u_0, \ldots, u_{n-1}$, and we can take spherical convex hull of $u_0, \ldots, u_{n-1}$, to be denoted $C$.
Again, since the $u_i$ are independent, the set $C$ has interior.  
Since $C$ is contained inside an open hemisphere, $C^*$ also has interior.  By the previous lemma, 
$C^* \cap -C$ is non-empty and open.  By the density of regular normals, there exists regular
unit normal $u$ for $K$ such that $u$ lies in the interior of  $C^* \cap -C$.   
Since $u$ lies in the interior of $C^*$, each $u \cdot u_i < 0$, so that
$u \notin N_i$ for any $i$ (by our choice of each $u_i \in N_i$).
It follows that $x = K^u$ is distinct from the previous exposed points $x_0, \ldots, x_{n-1}$.
Moreover, since $u$ lies in the interior of $-C$
$$-u = a_0 u_0 +\cdots + a_{n-1}u_{n-1}$$
for some $a_i > 0$, so that
$$a_0 u_0 +\cdots + a_{n-1}u_{n-1} + u = o.$$
Set $u_n = u$ and $x_n = x$.  The Minkowski existence theorem \cite[p. 125]{Bonn2}\cite[p. 390]{red} 
(or a much simpler Cramer's rule argument) yields 
an $n$-simplex with unit normals $u_0, \ldots, u_n$.  Scaling this simplex to
circumscribe $K$, each $i$th facet will meet the boundary of $K$ at exactly the distinct
exposed point $x_i$.  
\end{proof}

\begin{theorem} If $K \in \KK_n$ has at least $n+1$ exposed points, 
then there exists a simplex $S \in \KK_n$ such that 
each projection $S_u$ contains a translate of the projection $K_u$, 
while $S$ does not contain a translate of $K$.
\label{toomany}
\end{theorem}

\begin{proof}
If $\dim(K) = n$ then Theorem~\ref{toomany} immediately follows from Theorem~\ref{got-normals}
and Lemma~\ref{um}.

If $\dim(K) = d < n$, let $\xi$ denote the affine hull of $K$. 
By Theorem~\ref{got-normals}, there exists a
a $d$-dimensional simplex $Q \subseteq \xi$ 
that circumscribes $K$ in $\xi$ and whose $d+1$ facet unit normals are regular unit normals of $K$.
Since $K$ has $n+1$ exposed points, there are (at least) another $n-d$ regular unit normals of $K$
(in $\xi$) at these additional exposed points.
After intersecting $Q$ with supporting half-spaces (in $\xi$) of $K$ relative to these additional
$n-d$ normals, we obtain a polytope $Q_1$ in $\xi$ whose $n+1$ facet unit normals are regular unit normals of $K$.
Since $\dim \xi < n$, apply
small perturbations of these $n+1$ facet unit normals to $Q$
along $\xi^\perp$ to obtain facet normals of a simplex $S$ in $\RR^n$, whose facet normals are
still regular unit normals to $K$ in $\RR^n$.  

In either instance, we have obtained a simplex $S \supseteq K$, 
so that $K$ meets the boundary of $S$ at exactly $n+1$ points, one point from the relative interior 
of each facet of $S$.  By Lemma~\ref{um} 
there exists $\epsilon > 1$ such that $\epsilon K_u$ can be translated inside
$S_u$ for all $u$.  But $\epsilon K$ cannot be translated inside $S$, since $S$ circumscribes
$K$ already, and $\epsilon > 1$.
\end{proof}

\begin{corollary} If $K \in \KK_n$ and $\dim(K) = n$, then there exists $L \in \KK_n$ such that 
each projection $L_u$ contains a translate of the projection $K_u$, 
while $L$ does not contain a translate of $K$.
\label{onlyflat}
\end{corollary}

The following proposition addresses an ambiguity regarding when 
shadows cover inside a larger ambient space.
\begin{proposition} Suppose that $\xi$ is a linear flat in $\RR^n$.
Let $K$ and $L$ be compact convex sets in $\xi$.  
Suppose that, for each $d$-subspace $\eta \subseteq \xi$, the projection
$L_\eta$ contains a translate of $K_\eta$.  Then
$L_\eta$ contains a translate of $K_\eta$ for every $d$-subspace $\eta \subseteq \RR^n$.
\label{include}
\end{proposition}

\begin{proof}  
Suppose that $\eta$ is a $d$-subspace of $\RR^n$.  Let $\hat{\eta}$ denote the orthogonal projection
of $\eta$ into $\xi$.  Since $\dim(\hat{\eta}) \leq \dim(\eta) = d$, we can
translate $K$ and $L$ inside $\xi$ so that $K_{\hat{\eta}} \subseteq L_{\hat{\eta}}$.  Let us assume this
translation has taken place.  Note that, for $v \in \hat{\eta}$, we now have
$h_K(v) \leq h_L(v)$.

If $u \in \eta$, then express $u = u_{\xi} + u_{\xi^\perp}$.  Since $K \subseteq \xi$,
$$h_K(u) = \max_{x \in K} x \cdot u 
=\max_{x \in K} x \cdot u_\xi = h_K(u_\xi),
$$
and similarly for $L$.  But since $u \in \eta$, we have $u_\xi \in \hat{\eta}$, so that
$$h_K(u) = h_K(u_\xi)
\leq h_L(u_\xi) = h_L(u).
$$
In other words, $K_\eta \subseteq L_\eta$.
\end{proof}

Theorem~\ref{toomany} can now be generalized.
\begin{theorem} Suppose that $d \in \{1, 2, \ldots, n-1\}.$
If $K$ has at least $d+2$ exposed points, then
there exists $L \in \KK_n$ such that 
the projection $L_\xi$ contains a translate of the projection $K_\xi$ 
for each $d$-dimensional subspace $\xi$,
while $L$ does not contain a translate of $K$.
\label{toomany-d}
\end{theorem}

\begin{proof} Note that $n > d$.
If $n = d+1$ then Theorem~\ref{toomany} applies, and we are done.

Suppose that Theorem~\ref{toomany-d} holds when $n=d+i$ for some $i \geq 1$.  If $n = d+i+1$,
then there are two possible cases to consider.

First, if $\dim K = n$, then Corollary~\ref{onlyflat} yields $L \in \KK_n$ such that
every shadow $L_u$ contains a translate of $K_u$, while $L$ does not contain a translate of $K$.
Since every $d$-subspace $\xi$ is contained in some hyperplane $u^\perp$, it follows {\it a fortiori}
that every $d$-dimensional shadow $L_\xi$ contains a translate of $K_\xi$ as well.

Second, if $\dim K < n$, the induction hypothesis holds in the (lower dimensional) affine hull 
$\mathrm{Aff}(K)$ of $K$.  In other words, 
there exists a compact convex set $L$ in $\mathrm{Aff}(K)$ such that 
the projection $L_\xi$ contains a translate of the projection $K_\xi$ 
for each $d$-dimensional subspace $\xi$ of $\mathrm{Aff}(K)$,
while $L$ does not contain a translate of $K$.  Since $\mathrm{Aff}(K)$ is a flat in $\RR^n$,
inclusion of $L$ in $\RR^n$ preserves these covering
properties, by Proposition~\ref{include}.
\end{proof}

\begin{corollary} If $\dim K = d+1$, where $d \leq n-1$, then there exists $L \in \KK_n$ such that 
the projection $L_\xi$ contains a translate of the projection $K_\xi$ 
for each $d$-dimensional subspace $\xi$,
while $L$ does not contain a translate of $K$.
\label{only-d}
\end{corollary}

\begin{proof} If $\dim K = d+1$ then $K$ must have at least $d+2$ exposed points \cite[p. 89]{Webster}, 
so that Theorem~\ref{toomany-d} applies.
 \end{proof}

\section{Concluding remarks}

Although we have restricted our covering questions to shadows given by {\em orthogonal} projections,
the next proposition shows that
the same results will apply when more general (possibly oblique) linear projections are admitted.
\begin{proposition} Let $K, L \in \KK_n$.  Let $\psi: \RR^n \rightarrow \RR^n$ be a
nonsingular linear transformation.  Then $L_u$ contains a translate of $K_u$ for all unit directions $u$
if and only if $(\psi L)_u$ contains a translate of $(\psi K)_u$ for all $u$.
\end{proposition}

\begin{proof} 
For $S \subseteq \RR^n$ and a nonzero vector $u$, let $\LL_S(u)$ denote the set of straight lines in $\RR^n$
parallel to $u$ and meeting the set $S$.
The projection $L_u$ contains a translate $K_u$ for each unit vector $u$ if and only if,
for each $u$, there exists $v_u$ such that
\begin{equation}
\LL_{K+v_u} (u) \subseteq \LL_L(u).
\label{hum}
\end{equation}
But $\LL_{K+v_u}(u) = \LL_{K}(u) + v_u$ 
and $\psi \LL_{K}(u) = \LL_{\psi K}(\psi u)$.
It follows that~(\ref{hum}) holds if and only if
$\LL_{K}(u) + v_u \subseteq \LL_L(u)$,
which in turn holds if and only if
$$\LL_{\psi K}( \psi u) + \psi v_u \subseteq \LL_{\psi L}(\psi u) \;\;\; \hbox{ for all unit } u.$$
Set 
$$\tilde{u} = \frac{\psi u}{|\psi u|} \;\;\; \hbox{ and } \;\;\; \tilde{v} = \psi v_u.$$
The relation~(\ref{hum}) now holds if and only if,
for all $\tilde{u}$, there exists $\tilde{v}$ such that
$$\LL_{\psi K}(\tilde{u}) + \tilde{v} \subseteq \LL_{\psi L}(\tilde{u}),$$
which holds if and only if 
$(\psi L)_{\tilde{u}}$ contains a translate of $(\psi K)_{\tilde{u}}$ for all $\tilde{u}$.
\end{proof}

In this note we have addressed the existence of a compact convex set $L$, 
whose shadows can cover those of a given set $K$, without containing a translate of $K$ itself.
A reverse question is addressed in \cite{Klain-Circ}:  Given a body $L$, does there necessarily exist $K$ so that
the shadows of $L$ can cover those of $K$, while $L$ does not contain a translate of $K$?
A body $L$ is called $d$-{\em decomposable} if $L$ is a {\em direct} Minkowski sum (affine Cartesian
product) of two or more convex bodies each of dimension at most $d$.  A body 
$L$ is called $d$-{\em reliable} if, whenever each $d$-shadow of $K$ can be translated inside
the corresponding shadow of $L$, it follows that $K$ can itself be translated inside $L$.
In \cite{Klain-Circ} it is shown that $d$-decomposability implies $d$-reliability, although
the converse is (usually) false.  The results
in \cite{Klain-Circ,Klain-Shadow}, along with those of the present article, motivate the following
related open questions:
\begin{enumerate}
\item[I.]  Under what symmetry (or other) 
conditions on a compact convex set $L$ in $\RR^n$ is $d$-reliability
equivalent to $d$-decomposability, for $d > 2$? 
\end{enumerate}
In \cite{Klain-Circ} it is shown that 1-reliability is equivalent to 1-decomposability.
That is, only parallelotopes are 1-reliable.   It is also shown that a centrally symmetric
compact convex set is 2-reliable if and only if it is 2-decomposable.  However, this equivalence fails 
for bodies that are not centrally symmetric.

Denote the $n$-dimensional (Euclidean) volume of $L \in \KK_n$ by $V_n(L)$.  
\begin{enumerate}
\item[II.]  
Let $K, L \in \KK_n$ such that $V_n(L) > 0$, and let $1 \leq d \leq n-1$.
Suppose that the orthogonal projection 
$L_\xi$ contains a translate of the projection
$K_\xi$ for all $d$-subspaces $\xi$ of $\RR^n$.  \\

\noindent
What is the best upper bound for the ratio $\frac{V_n(K)}{V_n(L)}$? 
\end{enumerate}
In \cite{Klain-Shadow} it shown that $V_n(K)$ may exceed $V_n(L)$, although $V_n(K) \leq n V_n(L)$.
This crude bound can surely be improved.

\begin{enumerate}
\item[III.] Let $K, L \in \KK_n$, and let $1 \leq d \leq n-1$.
Suppose that, for each $d$-subspace $\xi$ of $\RR^n$, the orthogonal projection 
$K_\xi$ of $K$ can be moved inside $L_\xi$ by some {\em rigid motion}
(i.e.~a combination of translations, rotations, and reflections).  \\

\noindent
Under what simple (easy to state, easy to verify) additional conditions
does it follow that $K$ can be moved inside $L$ by a rigid motion?
\end{enumerate}
Because of the non-commutative nature of rigid motions 
(as compared to translations), covering via rigid motions may be more difficult 
to characterize than the case in which only translation is allowed.

\bibliographystyle{amsplain}

\end{document}